\documentclass{amsart}

\usepackage{xypic, epsf}
\usepackage{amssymb}
\usepackage{ifthen}

\newtheorem{theorem}{Theorem}[section]
\newtheorem*{maintheorem}{Theorem~\ref{theorem:5C}}
\newtheorem{proposition}[theorem]{Proposition}
\newtheorem{lemma}[theorem]{Lemma}

\newtheorem{corollary}[theorem]{Corollary}

\newtheorem*{theorem*}{Theorem}
\newtheorem*{lemma*}{Lemma}
\newtheorem*{proposition*}{Proposition}
\newtheorem*{claim*}{Claim}
\newtheorem*{corollary*}{Corollary}
\theoremstyle{definition}
\newtheorem*{definition}{Definition}
\newtheorem*{convention}{Convention}

\newcommand{\cross}{\times}
\newcommand{\abs}[1]{\left\lvert{#1}\right\rvert}
\newcommand{\fabs}[1]{|{#1}|}
\newcommand{\set}[1]{\left\{#1\right\}}

\newcommand{\half}{\frac{1}{2}}

\newcommand{\Z}{{\mathbb{Z}}}

\newcommand{\del}{\partial}

\newcommand{\C}{{\mathcal C}}
\newcommand{\setvert}{\;\vert\;}
\renewcommand{\phi}{\varphi}
\newcommand{\inner}[1]{\left<{#1}\right>}
\newcommand{\gromov}[1]{\left({#1}\right)}

\newcommand{\groupgen}[1]{\left<{#1}\right>}
\newcommand{\fgroupgen}[1]{\langle{#1}\rangle}

\newcommand{\putrefhere}[1]{[{\bf ??}]}

\newcommand{\eof}{\$}
\newcommand{\normal}{\lhd}
\newcommand{\action}{\backslash}
\newcommand{\nbhd}[2]{\overline{N_{#1}(#2)}}

\newenvironment{roman-enumerate}{\begin{enumerate}}{\end{enumerate}}

\newenvironment{tlemma}[1][true]{\ifthenelse{\equal{#1}{true}}{\begin{lemma}}{\begin{lemma}\label{lemma:#1}}}{\end{lemma}}
\newenvironment{tlemmacite}[2][true]{\ifthenelse{\equal{#1}{true}}{\begin{lemma}[#2]}{\begin{lemma}[#2]\label{lemma:#1}}}{\end{lemma}}
\newenvironment{tproposition}[1][true]{\ifthenelse{\equal{#1}{true}}{\begin{proposition}}{\begin{proposition}\label{proposition:#1}}}{\end{proposition}}
\newenvironment{tpropositioncite}[2][true]{\ifthenelse{\equal{#1}{true}}{\begin{proposition}[#2]}{\begin{proposition}[#2]\label{proposition:#1}}}{\end{proposition}}
\newenvironment{ttheorem}[1][true]{\ifthenelse{\equal{#1}{true}}{\begin{theorem}}{\begin{theorem}\label{theorem:#1}}}{\end{theorem}}
\newenvironment{ttheoremcite}[2][true]{\ifthenelse{\equal{#1}{true}}{\begin{theorem}[#2]}{\begin{theorem}[#2]\label{theorem:#1}}}{\end{theorem}}

\DeclareMathOperator{\Homeo}{Homeo}

\DeclareMathOperator{\id}{id}
\DeclareMathOperator{\caseif}{if}

\title{Quasiconvex Subgroups and Nets in Hyperbolic~Groups}
\author{Thomas Mack}
\address{Department of Mathematics \\Caltech\\ Pasadena, CA 91125}
\email{tmack@caltech.edu}

\begin{document}
\begin{abstract}
Consider a hyperbolic group $G$ and a quasiconvex subgroup $H\subset G$ with $[G:H] = \infty$.  We construct a set-theoretic section $s:G/H \to G$ of the quotient map (of sets) $G\to G/H$ such that $s(G/H)$ is a net in $G$; that is, any element of $G$ is a bounded distance from $s(G/H)$.  This section arises naturally as a set of points minimizing word-length in each fixed coset $gH$.  The left action of $G$ on $G/H$ induces an action on $s(G/H)$, which we use to prove that $H$ contains no infinite subgroups normal in $G$.
\end{abstract}
\maketitle


\section{Introduction}\noindent
Let $G$ be a group with finite generating set $\Sigma$.  The Cayley graph $\Gamma = C(G, \Sigma)$ is defined to be the graph with vertex set $G$ and edges connecting those $g, g'\in G$ with $g = g'\sigma$ for some $\sigma\in \Sigma$.  (We assume that $\Sigma$ is closed under inversion, so that this relation is symmetric.)  A hyperbolic group is one for which $\Gamma$ has the large-scale structure of a tree.  That is, geodesics are ``almost'' unique, in the sense that there exists a constant $C > 0$ such that any two geodesics $\gamma(t), \gamma'(t)$ between the same points satisfy $d(\gamma(t), \gamma'(t)) < C$ for all $t$.  

Let $H$ be a finitely generated subgroup of $G$, and choose a finite generating set $\Sigma'$ for it.  Assume without loss of generality that $\Sigma'\subset \Sigma$.  Then the Cayley graph $\Gamma' = C(H, \Sigma')$ naturally embeds in $\Gamma$, and so there are two natural metrics on $\Gamma'\subset \Gamma$: the path-length metric $d'$ considering chains that remain in $\Gamma'$ for all time, and the path-length metric $d$ condering all chains in $\Gamma$.  A quasiconvex subgroup is one for which these two metrics differ by no more than a constant multiplicative factor.  

The purpose of this paper is to prove the following theorem:
\begin{maintheorem}
Let $G$ be a hyperbolic group, and let $H\subset G$ be a quasiconvex subgroup.  If $[G:H] = \infty$, then there exists a (set-theoretic) section $s:G/H \to G$ of the quotient map $G \to G/H$ such that $s(G/H)$ is a net in $G$; that is, $\sup_{g\in G} d(g, s(G/H))$ is finite.
\end{maintheorem}\noindent
The argument depends on showing that for suitable $s$, there exists a finite automaton recognizing the language $L$ of points in $s(G/H)$.  By a geometric argument, any point in $\Gamma$ lies within a bounded distance of the prefix closure $\overline{L}$ of $L$.  Since $L$ is a regular language, it follows that any point in $\overline{L}$ is a bounded distance from a point in $L$, proving the theorem.  

The first two preliminary sections of the paper summarize general results in hyperbolic geometry that are used in the subsequent sections.  Section~2 is a broad overview of hyperbolic topological spaces and hyperbolic groups.  General references for this section include \cite{gromov}, which discusses hyperbolic spaces and groups in detail; \cite{metric}, which covers the large-scale geometry of general metric spaces; and \cite{blue}, which outlines many results in combinatorial group theory used throughout the paper.  Section~3 is an overview of finite automata and regular languages.  This machinery is useful not only for its direct use in the proof of the main theorem, but also because arbitrary hyperbolic groups have an automatic structure; the set of geodesics in $\Gamma$ can be recognized by a finite automaton.  The general material on finite automata can be found in \cite{words}, and \cite{cannon} contains specific applications to hyperbolic groups.  Section~4 is a brief summary of the problem of finding nets as in the main theorem for arbitrary groups, including a few examples and counterexamples.  Section~5 contains the proof of the main theorem and a result about normal subgroups embedded in quasiconvex subgroups that follows from it.

The material in this paper is a summary of some of the results in the author's Ph.D. thesis, supervised by Danny Calegari.
\section{Hyperbolic Spaces and Groups}\noindent
In this section, we recall some basic results and machinery of large-scale and hyperbolic geometry.  We assume the reader is familiar with the fundamental elements of the theory of coarse geometry, including quasigeodesics, quasi-isometries, and quasi-isometric embeddings.  We also assume familiarity with the basic results of hyperbolic geometry, such as $\delta$-hyperbolicity of metric spaces or groups, the Morse Lemma, and the ideal boundary of a hyperbolic space.  General references for the subject, including the material in this section, include \cite{bridson}, \cite{asymptotic}, and \cite{gromov}.

Let $(X, d)$ be a geodesic metric space.  For any $\epsilon \geq 0$ and $S\subset X$, define the closed $\epsilon$-neighborhood $\nbhd{\epsilon}{S}$ around $S$ by
\begin{align*}
\nbhd{\epsilon}{S} = \set{x\in X:\, \text{$d(x, p) \leq \epsilon$ for some $p\in S$}}.
\end{align*}
We denote a geodesic (not necessarily unique) with endpoints $p$ and $q$ by $[pq]$.  Call $X$ {\it $\delta$-hyperbolic} if it satisfies both of the following conditions:
\begin{roman-enumerate}
\item For any points $p, q, r\in X$ and any geodesics $[pq], [qr], [pr]$ in $X$, we have $[pr]\subset \nbhd{\delta}{[pq]\cup [qr]}$. 
\item For all $p, q, r, x\in X$, we have
\begin{align*}
\gromov{p. q}_x \geq \min\set{\gromov{p. r}_x, \gromov{q. r}_x} - \delta,
\end{align*}
where $\gromov{\cdot. \cdot}_\cdot$ denotes the {\it Gromov product}
\begin{align}\label{gromov-product}
\gromov{p.q}_x
	&= \half \left(d(p, x) + d(q, x) - d(p, q)\right).
\end{align}
\end{roman-enumerate}
Call $X$ {\it hyperbolic} if it is $\delta$-hyperbolic for some $\delta \geq 0$.  Note that if $X$ satisfies the first condition for a particular value of $\delta$, then it also satisfies the second for some $\delta'$, and vice versa.  Recall the following fundamental lemma of hyperbolic geometry, which states that quasigeodesics in a hyperbolic space remain uniformly close to geodesics:
\begin{tlemmacite}[1M]{{Morse Lemma, \cite[8.4.20]{metric}}}
Let $(X, d)$ be a $\delta$-hyperbolic geodesic space.  For any $K > 0$ and $\epsilon \geq 0$, there exists a constant $C > 0$, depending only on $\delta$, $K$, and $\epsilon$, such that any $(K, \epsilon)$-quasigeodesic segment with endpoints $p, q\in X$ lies in the $C$-neighborhood of a geodesic $[pq]$.
\end{tlemmacite}
To simplify notation, we adopt the following convention throughout:
\begin{convention}
All generating sets for groups are assumed to be finite and closed under inversion.
\end{convention}
For any group $G$ with generating set $\Sigma$, define the {\it Cayley graph} $\Gamma = C(G, \Sigma)$ to be the graph with vertex set $G$ and edges connecting those vertices $g, g'$ with $g' = gc$ for some $c\in \Sigma$.  We often implicitly identify $G$ with the set of vertices in $\Gamma = C(G, \Sigma)$, or even with $\Gamma$ itself.  The free group $F(\Sigma)$ with basis $\Sigma$ admits a homomorphism onto $G$ sending a word $[\sigma_1]\cdots [\sigma_n]$ to the corresponding element $\sigma_1 \cdots \sigma_n$ of $G$.  Denote this evaluation map by either $w \to \pi(w)$ or $w \to \overline{w}$.  To simplify notation, we often write $\overline{w}$ or $\pi(w)$ simply as $w$ when the intended meaning is clear.  For any $g\in G$, the {\it length} of $G$ with respect to $\Sigma$ is defined to be
\begin{align*}
\abs{g}_\Sigma = \min \set{\abs{w}:\, \alpha\in \pi^{-1}(w)},
\end{align*}
where $\abs{w}$ is the usual word length in $F(\Sigma)$; explicitly, $\abs{1} = 0$ and 
\begin{align*}
\abs{\sigma_1^{n_1} \dots \sigma_r^{n_r}} = \abs{n_1} + \cdots + \abs{n_r}
\end{align*}
for any nontrivial reduced word in $F(\Sigma)$.  Metrize $\Gamma$ by setting $d(x, y) = \fabs{x^{-1}y}_\Sigma$ for all vertices $x,y \in G$, and extend this metric over the edges to make $\Gamma$ into a geodesic length space.  In particular, the length of a chain $(x_0, \dots, x_n)$ in the graph $\Gamma$ is simply $n$.  Call $G$ {\it (word- or Gromov-)hyperbolic} if the Cayley graph $C(G, \Sigma)$ is a hyperbolic geodesic space.  Since any other generating set $\Sigma'$ for $G$ produces a Cayley graph $C(G, \Sigma')$ quasi-isometric to $C(G, \Sigma)$, word-hyperbolicity is independent of the choice of $\Sigma$ by the Morse Lemma. 

Throughout this section, let $G$ denote a hyperbolic group with generating set $\Sigma$ and Cayley graph $\Gamma = C(G, \Sigma)$.  Since $\Gamma$ is a hyperbolic geodesic space with only finitely many geodesics between any two points, we can extend the Morse Lemma~\ref{lemma:1M} to geodesic rays on $\Gamma$.
\begin{tlemma}[1F]
Let $G$ be a $\delta$-hyperbolic group, and set $\Gamma = C(G, \Sigma)$ for some generating set $\Sigma$ of $G$.  For any $K > 0$ and $\epsilon \geq 0$, there exists a constant $C > 0$, depending only on $\delta, K,$ and $\epsilon$, such that any $(K, \epsilon)$-quasigeodesic ray from a point $p$ lies within a distance $C$ of a geodesic ray from $p$.
\end{tlemma}
Let $X$ be a hyperbolic space, and fix $x_0\in X$.  Let $\Omega(X)$ denote the set of sequences of points in $X$ with $\gromov{p_i.p_j}_{x_0} \to \infty$ as $i, j\to\infty$.  Write $(p_i) \sim (p'_i)$ if $\liminf_{i, j\to\infty} \gromov{p_i.p'_i}_{x_0} = \infty.$  By \eqref{gromov-product}, $\sim$ is an equivalence relation on $\Omega(X)$. The set of equivalence classes $\Omega(X)/\!\!\!\sim$ is the {\it hyperbolic boundary} $\del X$.  For a hyperbolic group $G$, set $\Omega(G) = \Omega(\Gamma)$ and $\del G = \del \Gamma$.  For any hyperbolic spaces $X, Y$ and any quasi-isometry $f:X\to Y$, let $f_\infty:\del X \to \del Y$ denote the map $f(p_i) = (fp_i)$.  The Gromov product on $X$ extends across the boundary $\del X$, and the resulting function defines a metric on $\del X$.  For any quasi-isometry $f:X \to Y$, the map $f_\infty : \del X \to \del Y$ is continuous in this topology.  In particular, the left action $L_g(g') = gg'$ of $G$ on itself induces a homomorphism $G \to \Homeo(\del G)$, given by $g \to (L_g)_\infty$.  This action of $G$ on $\del G$ yields the following three useful propositions:
\begin{tpropositioncite}[1H]{{\cite[V.58]{blue}}}
Let $g\in G$.  If $g$ is not torsion, then the centralizer $C_G(g)$ of $g$ is virtually cyclic.  (That is, $C_G(g)$ is a finite extension of a cyclic group.)
\end{tpropositioncite}
\begin{tpropositioncite}[1I]{{\cite[V.58]{blue}}}
If $G$ is not elementary, then it contains a nonabelian free group.
\end{tpropositioncite}
\begin{tpropositioncite}[1J]{{\cite[V.56]{blue}}}
There are only finitely many conjugacy classes of $G$ that consist of torsion.
\end{tpropositioncite}
Recall that a subgroup $H\subset G$ is {\it $K$-quasiconvex} for some $K\geq 0$ if for all $p, q\in H$, any geodesic $[pq]$ in $\Gamma$ lies in $\nbhd{K}{H}$.  If the particular value of $K$ is irrelevant, we refer to $H$ as simply {\it quasiconvex}.  Any quasiconvex subgroup of a hyperbolic group $G$ is also hyperbolic; in particular, it is finitely generated.  It follows immediately from the Morse Lemma that for any generating set $\Sigma'$ of $H$, the inclusion $C(H, \Sigma') \to C(G, \Sigma\cup \Sigma')$ is a quasi-isometric embedding.

\section{Finite Automata}
In this section, we recall Cannon's result that hyperbolic groups are automatic.  We assume that the reader is familiar with the basic constructions of combinatorial group theory, such as finite automata and automatic structures.  Most of the defintions and results in this section are taken from \cite{words}, which also covers the material outlined here in more detail.

Recall that a generalized finite automaton (or finite automaton, or simply automaton) is a quintuple $M = (S, \Sigma, \mu, Y, s_0)$ satisfying the following properties:
\begin{roman-enumerate}
\item The {\it state set} $S$ is a finite set.
\item The {\it alphabet} $\Sigma$ is a finite set.
\item The {\it transition function} $\mu$ is a map $S \cross \Sigma^*\to 2^S$ such that the language $\set{w\in \Sigma^*:\, s'\in \mu(s, w)}$ over $\Sigma$ is regular for each fixed $s, s'\in S$.
\item The {\it set of accept states} $Y$ is a subset of $S$.
\item The {\it set of start states} $S_0$ is a subset of $S$. 
\end{roman-enumerate}
Let $L(M)\subset \Sigma^*$ denote the language recognized by $M$.  For any automaton $M$, the language $L(M)$ is regular.  The connection between hyperbolic groups and finite automata is that words of minimum length (with respect to some fixed generating set) in a hyperbolic group form a regular language.  More precisely, consider an arbitrary finitely generated group with generating set $\Sigma$.  Denote the semigroup homomorphism $\Sigma^* \to X$ sending $[\sigma_1] \cdots [\sigma_n]$ to the corresponding element $\sigma_1 \cdots \sigma_n$ of $X$ by $w \to \pi(w)$ or $w \to \overline{w}$.  Call $X$ {\it automatic} if there exist finite automata $W$ over $\Sigma$ and $M_x$ over $\Sigma^2$ for each $x\in \Sigma\cup \set{\epsilon}$ that satisfy the following two properties:
\begin{roman-enumerate}
\item For each $g\in X$, there exists some $w\in L(W)$ with $\overline{w} = g$. 
\item The language $L(M_x) = \set{(u, v)\in L(W)\cross L(W):\, \overline{ux} = \overline{v}}$.
\end{roman-enumerate}
The automaton $W$ is called a {\it word acceptor}, and $M_x$ is called a {\it multiplier automaton} for $x \not = \epsilon$ and an {\it equality recognizer} for $x = \epsilon$.  The condition of having an automatic structure is called {\it automation}.  Automation is independent of the choice of generating set $\Sigma$.  Fix a hyperbolic group $G$ with generating set $\Sigma$ and Cayley graph $\Gamma = C(G, \Sigma)$.  The automatic structure on $G$ arises naturally from Cannon's idea of {\it cone types}~\cite{cannon}, defined below.
\begin{definition}
For any $x, y\in G$, write $x\leq y$ if there exists a geodesic segment $[1y]$ passing through $x$.  Define an equivalence relation on $G$ by setting $x\sim x'$ if the inequality $x\leq y$ holds exactly when $x'\leq (x'x^{-1})y.$  The quotient $\C(G) = G/\!\!\sim$ is the {\it set of cone types} of $G$, and the image of $x\in G$ in $\C(G)$ is the {\it cone type} of $x$, denoted by $C(x)$.
\end{definition}
\begin{ttheoremcite}[3F]{Cannon}
For any hyperbolic group $G$, the set of cone types $\C(G)$ is finite.
\end{ttheoremcite}
\begin{proof}
See \cite{cannon} for a geometric proof or \cite[3.2]{words} for a more combinatorial one.
\end{proof}
To simplify notation, write $u\leq v$ for words $u, v\in \Sigma^*$ if $\overline{u}\leq \overline{v}$.  For any $C\in \C(G)$, write $C\leq x$ if $g\leq gx$ for some (and hence every) $g\in G$ with $C(g) = C$.  For any fixed $x\in G$, the relation $g\leq gx$ depends only on the cone type of $g$; furthermore, if it holds, then $C(gx)$ depends only on $C(g)$.  Thus write $C(g)x$ for $C(gx)$ if $C(g)\leq x$.  Define the {\it language of geodesics} $\Lambda(G)$ to be the set of all words $w\in \Sigma^*$ with $\abs{w} = \abs{\overline{w}}_\Sigma.$  (Note that $\Lambda$ also depends on the choice of $\Sigma$.) Thus elements of $\Lambda(G)$ correspond to geodesics in $\Gamma$ from $1$ to any other point.  The language $\Lambda(G)$ over $\Sigma$ is regular for any hyperbolic group $G$.  Cannon's result~\ref{theorem:3F} therefore implies that any hyperbolic group has an automatic structure.
\begin{tpropositioncite}[3E]{{\cite[2.3.4, 3.2]{words}}}
Every $\delta$-hyperbolic group $G$ has an automatic structure with $L(W) = \Lambda(G)$.
\end{tpropositioncite}
Applying the Pumping Lemma to the automaton described in Proposition~\ref{proposition:3E} produces the following result:
\begin{tproposition}[3G]
If $G$ is infinite, then there exists a constant $C > 0$ such that for every $g\in G$, some geodesic ray $r$ from $1$ in $\Gamma$ passes through $\nbhd{C}{g}$.  
\end{tproposition}


\section{Nets in Groups}
Let $M$ be a metric space (not necessarily hyperbolic).  For any $C > 0$, a subspace $M'\subset M$ is a {\it $C$-net} if $M\subset \nbhd{C}{M'}$.  Call $M'$ a {\it net} in $M$ if $M'$ is a $C$-net for some $C$; or, equivalently, if $d(p, M')$ is bounded for all $p\in M$.  Similarly, for a finitely generated group $X$, a subgroup $X'\subset X$ is a {\it net} if it is a net in the metric space $C(X, \Sigma)$ for some finite generating set of $\Sigma$.  Note that we do not require that $X$ be hyperbolic or that $X'$ be finitely generated.  Since the Cayley graph $C(X, \Sigma)$ is independent of $\Sigma$ up to quasi-isometry, the condition of being a net in $X$ is independent of the particular choice of generating set $\Sigma$ for $X$.  Consider the problem of finding pairs $(X, X')$ with $X'\subset X$ that satisfy the following property:
\begin{align}\tag{$*$}\label{prop}
\text{\it There exists a section $s:X/X'\to X$ such that $s(X/X')$ is a net in $G$.}
\end{align}
In \eqref{prop}, $X/X'$ denotes the space of right cosets of $X'$ in $X$; we do not require $X'\normal X$.  In particular, the desired map $s:X/X' \to X$ is only a map of sets, not a group homomorphism or a continuous map. The goal of this paper is to prove that the pair $(G, H)$ satisfies \eqref{prop} if $G$ is hyperbolic and $H\subset G$ is a quasiconvex subgroup of infinite index.  In order to motivate this result, we provide a few examples and counterexamples of pairs satisfying \eqref{prop} in this section.

\begin{tlemma}[4A]
Let $1 \to N \to E \xrightarrow{\pi} Q \to 1$ be an exact sequence of groups.  Suppose $E$ is finitely generated.  Then the pair $(E, N)$ satisfies \eqref{prop} iff $N$ is finite.
\end{tlemma}
\begin{proof}
It is clear that $(E, N)$ satisfies \eqref{prop} if $N$ is finite.  Thus assume $N$ is infinite.  Fix generating sets $\Sigma$ and $\Sigma' = \pi(\Sigma)$ for $E$ and $Q$, respectively.  Suppose instead that there exists a section $s:Q\to E$ of $\pi$ such that $s(Q)$ is a $C$-net (with respect to $\Sigma$) in $C(E, \Sigma)$ for some constant $C$.  Clearly $\abs{\pi g}_{\Sigma'} \leq \abs{g}_{\Sigma}$ for all $g$.  Thus for any $g, g'\in E$, we have 
\begin{align}\label{pi}
d(g, g')
	&= \abs{g^{-1}g'}_{\Sigma}
	\geq \abs{(\pi g)^{-1} (\pi g')}_{\Sigma'}
	= d(\pi g, \pi g').	
\end{align}
Fix $q_0\in Q$.  Since $N$ is infinite, there exists some $g\in E$ with $g\in s(q_0) N$ and $d(g, s(q)) > C$ for all $q\in Q$ with $d(q, q_0) \leq C$.  By \eqref{pi}, $d(g, s(q)) > C$ for all $q\in Q$.  The lemma follows from this contradiction.
\end{proof}
Thus pairs $(E, N)$ with $N\normal E$ cannot satisfy \eqref{prop} because the preimages $\pi^{-1}(q)$ remain uniformly separated: $d(\pi^{-1}q, \pi^{-1}q') \geq d(q, q')$ for all $q, q'\in Q$.  Hence to construct a section $s:X/X' \to X$ of the required type, we need to consider subgroups $X'$ for which the cosets $gX'$ exhibit more complicated behavior.  It was proved in Lemma~\ref{lemma:4A} that for any exact sequence $1\to N \to E \to Q\to 1$, the pair $(E, N)$ does not satisfy \eqref{prop}.  If this sequence splits, then we can embed $Q$ in $E$ and consider the pair $(E, Q)$; that is, we consider $(E, Q)$ for $E$ a semidirect product $N\rtimes Q$.  In order to analyze this problem, it is useful to consider the projection $f:E \to s(Q)$ instead of $s$ itself.  Note that $f$ moves points a bounded distance iff $s(Q)$ is a net in $E$.  The following lemma makes this observation more precise.
\begin{tlemma}[4B]
Let $X = N \rtimes X'$ for some action of $X'$ on $N$, and let $\pi:X \to X/N = X'$ denote the quotient map.  Suppose $X$ is finitely generated.  Let $d$ and $d'$ denote the metrices on $C(X, \Sigma)$ and $C(X', \Sigma')$, respectively, for some fixed generating sets $\Sigma$ of $X$ and $\Sigma = \pi(\Sigma)$ of $X'$.  Then $(X, X')$ satisfies \eqref{prop} iff there exists a function $\phi:X \to N$ and a constant $C$ that satisfy the following two properties:
\begin{roman-enumerate}
\item If $\phi(\alpha) = \phi(\alpha')$ for $\alpha\not = \alpha'$ in $X$, then $d'(\pi \alpha, \pi \alpha') \leq C$.
\item The distance $d(\alpha, (\phi \alpha, \pi \alpha))\leq C$ for all $\alpha\in X$.
\end{roman-enumerate}
\end{tlemma}
\begin{proof}
Suppose first that $(X, X')$ satisfies \eqref{prop}.  Then there exists a section $s:N\to X$ such that $S = s(N)$ is a $\half C$-net in the Cayley graph $C(G)$ for some $C$.  For each $\alpha\in X$, choose some $f(\alpha)\in S$ with $d(\alpha, f(\alpha)) \leq C/2$.  Write $f(\alpha) = (\phi(\alpha), \psi(\alpha))$ for some functions $\phi:X \to X'$ and $\psi:X \to N$.  Since $s$ is a section, we have $\psi(\alpha) = \psi(\alpha')$ for some $\alpha, \alpha'\in X$ iff $\phi(\alpha) = \phi(\alpha')$.  But
\begin{align*}
d'(\psi \alpha, \pi \alpha)
	&= d(\pi f\alpha, \pi \alpha)
	\leq d(f\alpha, \alpha)
	\leq C/2,
\end{align*}
so $\phi$ satisfies property (i).

Since $\Sigma' = \pi(\Sigma)$, we have $\abs{\pi\alpha}_{\Sigma'} \leq \abs{\alpha}_{\Sigma}$ for all $\Sigma$.  Assume without loss of generality that $(1, \sigma')\in \Sigma$ for all $\sigma'\in \Sigma'$.  Then $\abs{(1, x)}_{\Sigma} \leq \abs{x}_{\Sigma'}$ for all $x\in X$.  It follows that $\abs{(1, x)}_{\Sigma} = \abs{x}_{\Sigma'}$ for all $x\in X$.  We therefore have
\begin{align}
d((n, x), (n, x'))
	&= \abs{(x^{-1}.n^{-1}, x^{-1})(n, x')}_{\Sigma} \nonumber\\
	&= \abs{(1, x^{-1}x')}_{\Sigma} \nonumber\\
	&=\abs{x^{-1}x'}_{\Sigma'} \nonumber\\
	&= d'(x, x')\label{eq:projection}
\end{align}
for all $n\in N$ and $x\in X'$.  Hence for all $\alpha\in X$,
\begin{align*}
d(\alpha, (\phi \alpha, \pi \alpha))
	&\leq d(\alpha, (\phi \alpha, \psi \alpha)) + d_\Sigma((\phi \alpha, \psi \alpha), (\phi \alpha, \pi \alpha)) \\
	&= d'(\psi \alpha, \pi \alpha) + C/2 \\
	&= d'(\pi f\alpha, \pi \alpha) + C/2 \\
	&\leq d(f\alpha, \alpha) + C/2 \\
	&\leq C.
\end{align*}
The conclusion of the lemma therefore holds.

Conversely, suppose that such a function $\phi$ exists. Consider the set 
\begin{align*}
S = \set{(\phi(n, x), x):\, (n, x)\in X}.
\end{align*}
By property (ii), $S$ is a $C$-net in $\Gamma = C(G, \Sigma)$.  For any points $(n, x), (n, x')\in X$ for fixed $n\in N$, we have $d'(x, x') \leq C$ by (i).  Hence $d((n, x), (n, x')) \leq C$ by \eqref{eq:projection}.  Let $S'$ be a subset of $S$ such that $S'\cap (n\cross X)$ contains at most one point for each $n\in N$. Then $S'$ is a $2C$-net in $\Gamma$. Choose an arbitrary (set-theoretic) section $s:N\to X$ such that $S'\cap (n\cross X) = \set{s(n)}$ whenever this intersection is nonempty.  Then $s(N)\supset S'$ is also a $2C$-net in $\Gamma$, so $(X, X')$ satisfies \eqref{prop}.
\end{proof}
Using Lemma~\ref{lemma:4B}, we now describe a method for using sections to satisfying \eqref{prop} to construct such sections over larger groups.
\begin{tlemma}
Let $Q$ be a finitely generated group, and let $Q$ act on groups $N$ and $N'$.  Suppose $E = N\rtimes Q$ and $E' = N'\rtimes Q$ are finitely generated.  Set $X = (N\cross N')\rtimes Q$, where $Q$ acts on $N\cross N'$ via the diagonal map.  If $(E, Q)$ and $(E', Q')$ satisfy \eqref{prop}, then $(X, Q)$ also satisfies it.
\end{tlemma}  
\begin{proof}
Choose finite generating sets $\Sigma$ and $\Sigma'$ for $E$ and $E'$, respectively.  Assume without loss of generality that $\Sigma$ contains $(1, q)$ for any $(n, q)\in \Sigma$, and similarly for $\Sigma'$.  Let $(n\cross n', q)\in X$, and let $(n, q) = (n_1, q_1) \cdots (n_r, q_r)$ and $(n', q) = (n'_1, q'_1), \dots, (n'_s, q'_s)$ with all $(n_i, q_i)\in E$ and $(n'_i, q'_i)\in E'$.  Then
\begin{align*}
(n\cross n', q)
	&= (1\cross n_1, q_1) \cdots (n_r, q_r)(1, q_r^{-1}) \cdots (1, q_1^{-1})(n'_1\cross 1, q'_1) \cdots (n'_s \cross 1, q'_s).
\end{align*}
Hence $X$ is generated by 
\begin{align*}
S =& \set{(n\cross 1, q):\, (n, q)\in \Sigma} \cup \set{(1\cross n, q):\, (n, q)\in \Sigma'};
\end{align*}
furthermore, we have 
\begin{align*}
\abs{(n\cross n', q)}_{\Sigma\cross \Sigma'} \leq 2{(n, q)}_{\Sigma} + \abs{(n', q)}_{\Sigma'}.
\end{align*}
Equivalently, 
\begin{align}\label{eq:bound}
\abs{\alpha}_{\Sigma\cross \Sigma'} \leq 2\abs{\rho \alpha}_{\Sigma} + \abs{\rho' \alpha}_{\Sigma'}
\end{align}
for all $\alpha\in X$, where $\rho$ and $\rho'$ denote the quotient maps $\rho:X \to X/N' = E$ and $\rho':X \to X/N = E'$.

Suppose $(E, Q)$ and $(E', Q)$ satisfy \eqref{prop}.  Then there exist functions $\phi:E \to N$ and $\phi':E'\to N'$ satisfying properties (i) and (ii) in Lemma~\ref{lemma:4B} for some constant $C$.  Consider the function $\psi:X \to N\cross N'$ defined by $\psi = \phi \rho \cross \phi' \rho'.$  To simplify notation, denote the three quotient maps $E \to E/N = Q, E' \to E'/N' = Q,$ and $X \to X/(N\cross N') = Q$ by $\pi$.  If $\psi(\alpha) = \psi(\beta)$, then property (i) forces
\begin{align}\label{eq:prop-i}
d_{\pi(S)}(\pi \alpha, \pi \beta) \leq d_{\pi(\Sigma)}(\pi \rho \alpha, \pi \rho \beta) \leq C.
\end{align}
By property (ii), any $\alpha\in X$ satisfies
\begin{align}
d_S(\alpha, (\psi \alpha, \pi \alpha))
	&\leq 2 d_{\Sigma}(\rho \alpha, (\rho \psi \alpha, \pi \alpha)) + d_{\Sigma'}(\rho' \alpha, (\rho' \psi \alpha, \pi \alpha))\nonumber \\
	&\leq 2 d_{\Sigma}(\rho \alpha, (\phi \rho \alpha, \pi \rho \alpha)) + d_{\Sigma'}(\rho' \alpha, (\phi \rho' \alpha, \pi \rho' \alpha))\nonumber \\
	&\leq 3C\label{eq:prop-ii}
\end{align}
by \eqref{eq:bound}.  Combining \eqref{eq:prop-i} and \eqref{eq:prop-ii} shows that $(X, Q)$ satisfies \eqref{prop} by Lemma~\ref{lemma:4B}.
\end{proof}


\section{Nets in Hyperbolic Groups}\noindent
In this section, we prove that for any hyperbolic group $G$ and any quasiconvex subgroup $H\subset G$ with $[G:H]$ infinite, the pair $(G, H)$ satisfies \eqref{prop}; that is, there exists a section $s:G/H \to G$ of the quotient map (of sets) $G\to G/H$ such that $s(G/H)$ forms a net in the Cayley graph of $G$.  Fix a $\delta$-hyperbolic group $G$, a generating set $\Sigma$ of $G$, and a $K$-quasiconvex subgroup $H$.  Let $\Gamma = C(G, \Sigma)$, and let $d(x, y) = \fabs{x^{-1}y}_\Sigma$ denote the metric on $\Gamma$.  Set $\Lambda = \Lambda(G)$ and $\C = \C(G)$.  To simplify notation, abbreviate $\abs{\cdot}_\Sigma$ as $\abs{\cdot}$.  

\begin{figure}[t]
\begin{center}
\epsfbox{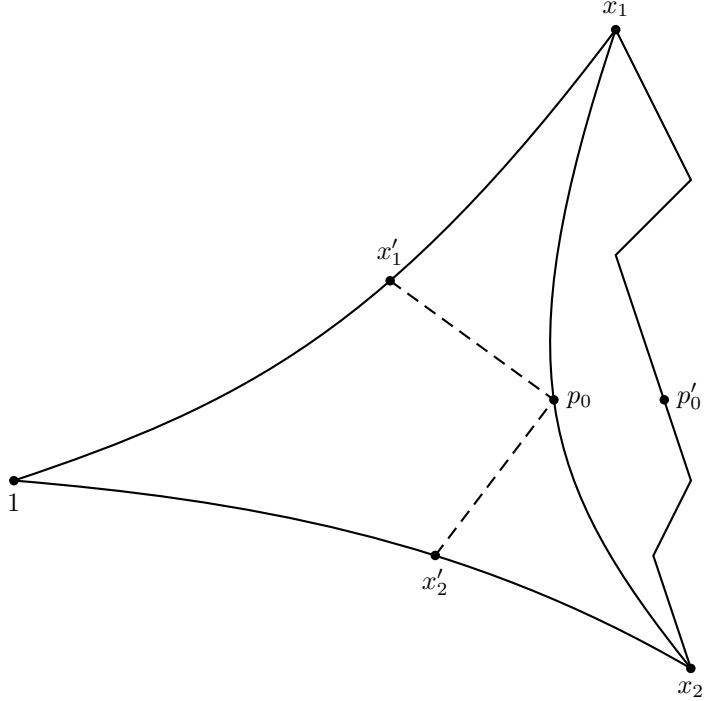}
\caption{The construction in Lemma~\ref{lemma:5A}.  The dotted lines each have length at most $\delta + 2$, and the quasigeodesic containing $p'_0$ stays in the $K$-neighborhood of the geodesic containing $p_0$.}\label{figure:lemmadiagram}
\end{center}
\end{figure}

For each $g\in G$, set $\sigma(g) = \min \set{\abs{gh}:\, h\in H}$ and $S = \set{g\in G:\, \abs{g} = \sigma(g)}.$  We form the desired section simply by choosing one point in each coset $S\cap gH$.  The diameter of the set $S\cap gH$ is uniformly bounded for all $g\in G$, so it suffices to prove that $S$ itself is a net. The crucial step is doing so is showing that the set of geodesics to points in $S$ is a regular language.  The prefix closure $\overline{S}$ of $S$ is then a finite distance in $\Gamma$ away from $S$ itself.  Thus if $\overline{S}$ is a net in $\Gamma$, then so is $S$.  As a subset of $\Gamma$, the set $\overline{S}$ consists of all points that lie on geodesic rays from $1$ that intersect $S$.  We use the quasiconvexity of $H$ to prove that any point in $\Gamma$ is a bounded distance from such a geodesic if $[G:H]=\infty$, completing the proof.  The first step in this argument is the following lemma, which provides a convenient bound or estimate for the distance between points on the same coset $gH$.  
\begin{tlemma}[5A]
Let $H\subset G$ be a quasiconvex subgroup, and fix some $g\in G$.  For all $x_1, x_2\in gH$, we have $d(x_1, x_2)\leq C_1 + \abs{x_1} + \abs{x_2} - 2 \sigma(g)$, where $C_1 > 0$ is a constant depending only on $G$ and $H$.
\end{tlemma}
\begin{proof}
Choose geodesics $[1x_1], [1x_2]$, and $[x_1x_2]$.  Set $f(p) = d(p, [1x_1]) - d(p, [1x_2])$ for all $p\in [x_1 x_2].$  For any two adjacent vertices $p$ and  $p'$, we have $\abs{f(p) - f(p')} \leq 2$.  Since $f(x_1) \leq 0$ and $f(x_2) \geq 0,$ it follows that there exists some $p_0\in [x_1x_2]$ with $\abs{d(p_0, [1x_1]) - d(p_0, [1x_2])} \leq 2.$  By the $\delta$-hyperbolicity of $\Gamma$, we have 
\begin{align*}
\min \set{d(p, [1x_1]), d(p, [1x_2])}\leq \delta
\end{align*}
for each $p\in [x_1 x_2].$  Thus $d(p_0, [1x_1]), d(p_0, [1x_2])\leq \delta + 2$; choose $x'_1\in [1x_1]$ and $x'_2\in [1x_2]$ realizing these inequalities.  The quasiconvexity of $H$ implies that there exists a point $p'_0\in gH$ with $d(p_0, p'_0) < K$.  See Figure~\ref{figure:lemmadiagram} for an illustration of this construction.  We have
\begin{align*}
\sigma(g)
	&\leq \abs{p'_0} \\
	&\leq \abs{x'_i} + d(x'_i, p_0) + d(p_0, p'_0) \\
	&\leq \abs{x'_i} + \delta + K + 2 \\
	&= \abs{x_i} - d(x_i, x'_i) + \delta + K + 2.
\end{align*}
Thus $d(x_i, x'_i) \leq \abs{x_i} - \sigma(g) + \delta + K + 2.$  We therefore have
\begin{align*}
d(x_1, x_2)
	&\leq d(x_1, x'_1) + d(x'_1, p_0) + d(p_0, x'_2) + d(x'_2, x_2) \\
	&\leq 2\delta + 4 + d(x_1, x'_1) + d(x_2, x'_2) \\
	&\leq 4\delta + 2K + 8 + \abs{x_1} + \abs{x_2} - 2\sigma(g),
\end{align*}
as required.
\end{proof}
For $x_1, x_2\in S$, Lemma~\ref{lemma:5A} forces $d(x_1, x_2)\leq C_1$.  The intersections $S\cap gH$ for $g\in G$ thus have bounded diameter.  In defining the section $s$ by choosing one point in $S\cap gH$ for each coset $gH$, the particular choice of points is therefore irrelevant in the large-scale geometry of $\Gamma$.  In particular, the condition that $s$ be a net is independent of this choice.

Set $\inner{g, g'} = 2(g^{-1}.g') = \abs{g} + \abs{g'} - \abs{gg'}$ for $g, g'\in G.$  Following the conventions in previous sections, we also write $\inner{x, y}$ for $\inner{\overline{x}, \overline{y}}$ with $x, y\in \Sigma^*$ to simplify notation.  The set $S$ consists precisely of those points $g\in G$ with $\abs{g} \leq \abs{gh}$ for all $h$.  Thus
\begin{align*}
S &= \set{g\in G:\, \text{$\inner{g, gh} \leq \abs{h}$ for all $h\in H$}}.
\end{align*}
In order to prove that $S$ is regular, we first prove that the condition $\inner{x, y} = n$ defines a regular language $L_n\subset \Lambda^2$ for each $n$.  The idea of the proof is to split each $y$ with $(x, y)\in L_n$ into subwords $w, c, w'\in \Sigma^*$ such that $\inner{x, y} = \inner{x, wcw'} = \inner{xy, c} + \inner{xwc, w'}$ with $\inner{x, wcw'}, \inner{xy, c} < n.$ By showing that the set of such words $(y, w, c, w')$ is regular, we therefore conclude that $L_n$ is regular by induction on $n$.
\begin{tlemma}[5B]
For any fixed $n\geq 0$, the language $L_n = \set{(x, y)\in \Lambda^2:\, \inner{x, y} = n}$ over $\Sigma^2$ is regular.
\end{tlemma}
\begin{proof}
Since $G$ is hyperbolic, $\Lambda$ is a regular language over $\Sigma$.  Hence there exists a determinstic finite automaton $M$ over $\Sigma$ with $L(M) = \Lambda$.  Let $M'$ denote the automaton obtained by replacing each arrow labelled $c\in \Sigma$ by arrows labelled $(c, \epsilon)$ and $(\eof, c)$ (with the same source and target), where $\eof$ denotes a padding character.  The language recognized by $M'$ consists of all $x, y\in \Sigma^*$ such that $xy\in \Lambda.$  Thus $L_0 = L(M') \cap \Lambda^2$, and so is regular.  Fix some positive integer $n$, and suppose that $L_m$ is regular for all $m < n.$  For any $x\in \Lambda$ and $c\in \Sigma,$
\begin{align*}
\inner{x, c} &= 
	\abs{x} + 1 - \abs{xc} 
	= \begin{cases}
		0 & \caseif x\leq xc; \\
		2 & \caseif x\geq xc; \\
		1 & \text{otherwise}.
	\end{cases}
\end{align*}
Set $P_i(c) = \set{x\in \Lambda\setvert \inner{x, c} = i}$ for $i = 0, 1, 2$ and $c\in \Sigma$.  Let $\C_c$ denote the set of cones $C$ such that $x\leq xc$ if $C(x) = C$.  Then
\begin{align*}
P_0(c) &= \set{x\in \Lambda\setvert C(x)\in \C_c} \\
P_2(c) &= \set{x\in \Lambda\setvert \exists y\in \Lambda:\, \overline{xc} = \overline{y}, C(y) \in \C_{c^{-1}}} \\
P_1(c) &= \Lambda \cap \neg(L_0 \cup L_2)
\end{align*}
Since $\C_c$ is finite, all three languages $P_i$ are regular.

Let $R\subset \Lambda^4$ consist of all quadruples $(x, y, z, w)$ with $y = y_1\cdots y_n$, $z = xy_1 \cdots y_i$, and $w = y_{i+1} \cdots y_n$, where $i$ is the largest index such that $xy_1 \cdots y_i\in \Lambda$.  We claim that the language $R$ is regular.  Construct a finite automaton $M$ over $\Sigma^4$ as follows.  Let $M$ have one state $s_C$ for each cone type $C\in \C$, and add two states $s_r$ and $s_a$.  Designate $s_a$ as the only accept state.  For each state $s_C$, add arrows labelled $(u, \epsilon, u, \epsilon)$ and $(\eof, u, u, \epsilon)$ from $s_C$ to the state $s_{Cu}$ for each $u\in \Sigma$ with $x\leq xu$ for $C(x) = C.$  For any other $u\in \Sigma$, attach an arrow from $s_C$ to $s_r$ labelled $(\eof, u, \epsilon, u).$  Add an arrow from $s_r$ to itself labelled $(\eof, u, \eof, u)$ for each $u\in \Sigma.$  Attach an arrow labelled $(\eof, \eof, \eof, \eof)$ from every state in $M$ to $s_a$.  It is clear that for any geodesics $x, y, z, w\in \Lambda$, the resulting automaton accepts the quadruple $(x, y, z, w)$ iff it lies in $R$.  Thus $R = L(M') \cap \Lambda^4$.  The claim follows.  

Let $R_c$ denote the set of quadruples $(x, y, z, w)\in R$ with $w_1 = c$, and let $R'_c$ denote the set of sextuples $(x, y, z, w, z', w')\in \Lambda^6$ such that $(x, y, z, w)\in R_c, \overline{z'} = \overline{zc}$, and $w' = w_2 \cdots w_{\abs{w}}.$  Since $R$ is regular, it is clear that $R_c$ and $R'_c$ are also regular for each $c\in \Sigma.$  For any $x, y, z\in G$, we have
\begin{align}
\inner{x, y} + \inner{xy, z} - \inner{y, z}
	&= \left(\abs{x} + \abs{y} - \abs{xy}\right) + \left(\abs{xy} + \abs{z} - \abs{xyz}\right) - \left(\abs{y} + \abs{z} - \abs{yz}\right) \nonumber\\
	&= \abs{x} + \abs{yz} - \abs{xyz} \nonumber\\
	&= \inner{x, yz}.\label{eq:induction}
\end{align}
Fix $x, y\in \Lambda$, and write $y = y_1 \cdots y_n$ with each $y_i\in \Sigma.$  Let $k\geq 0$ denote the largest index such that $xy_1 \cdots y_k\in \Lambda$.  By \eqref{eq:induction},
\begin{align*}
\inner{x, y} 
	&= \inner{xy_1\cdots y_{k+1}, y_{k+2} \cdots y_n} + \sum_{i=0}^k \inner{xy_1 \cdots y_i, y_{i+1}} \\
	&= \inner{xy_1\cdots y_k, y_{k+1}} + \inner{xy_1 \cdots y_{k+1}, y_{k+2}\cdots y_n}.
\end{align*}
Since $xy_1 \cdots y_{k+1}\not\in \Lambda,$ we have $\inner{xy_1 \cdots y_k, y_{k+1}} > 0$.  Thus $\inner{x, y} = n$ iff we have 
\begin{align*}
\inner{xy_1 \cdots y_k, y_{k+1}} &= i; & \inner{xy_1 \cdots y_{k+1}, y_{k+2}\cdots y_n} &= n - i
\end{align*}
for $i = 1$ or $2$.  Hence $L_n$ satisfies
\begin{align*}
L_n 
	= \bigcup_{c\in \Sigma; i = 1, 2}\Big\{(x, y)\in \Lambda^2&\setvert \exists (z, w, z', w')\in \Lambda^4: \\
	&(x, y, z, w, z', w')\in R'_c, \inner{z, c} = i, \inner{z', w'} = n - i\Big\}; \\
\intertext{that is,}
L_n
	= \bigcup_{c\in \Sigma; i = 1, 2}\Big\{(x, y)\in \Lambda^2&\setvert \exists (z, w, z', w')\in \Lambda^4: \\
	&(x, y, z, w, z', w')\in R'_c, z\in P_i(c), (z', w')\in L_{n-i}\Big\},
\end{align*}
Thus $L_n$ is regular.  The lemma follows by induction on $n$.
\end{proof}
\begin{tlemma}[5I]
The language $L = \set{x\in \Lambda:\, \overline{x}\in S}$ is regular.
\end{tlemma}
\begin{proof}
Assume without loss of generality that the generating set $\Sigma$ of $G$ contains a generating set $\Sigma'$ for the hyperbolic group $H$.  The language $\Lambda'\subset \Sigma'^* \subset \Sigma^*$ of geodesics in $C(H, \Sigma')$ is regular by the fact that $H$ is hyperbolic.  For any $g\in S$ and $h\in H$, we have 
\begin{align*}
\inner{g, gh} 
	&= \abs{g} + d(g, gh) - \abs{gh}
	\leq \abs{g} + \left(C_1 + \abs{gh} - \abs{g}\right) - \abs{gh}
	= C_1
\end{align*}
for some constant $C_1$ by Lemma~\ref{lemma:5A}.  Hence
\begin{align*}
S 
	&= \set{x\in \Lambda\setvert\forall y\in \Lambda':\, \inner{x, y} \leq \abs{y}}
	= \set{x\in \Lambda\setvert\forall y\in \Lambda':\, \inner{x, y}\leq \min(\abs{y}, C_1)}.
\end{align*}
Thus
\begin{align}
S
	&= \set{x\in \Lambda\setvert\forall y\in \Lambda':\, \inner{x, y}\leq C_1} \cap 
		\bigcup_{r=0}^{C_1} \left(\bigcap_{y\in \Lambda' \cap \nbhd{r}{1}} \set{x\in \Lambda:\, \inner{x, y}\leq r}\right)\label{eq:regularity-of-S}.
\end{align}
By Lemma~\ref{lemma:5B}, the language 
\begin{align*}
\set{(x, y)\in \Lambda\cross \Lambda':\, \inner{x, y}\leq n} = (\Lambda\cross\Lambda') \cap \bigcup_{i=0}^n \set{(x, y)\in \Lambda^2:\, \inner{x, y} = i}
\end{align*}
is regular for all $n.$  Thus $S$ is regular by \eqref{eq:regularity-of-S}.
\end{proof}
We now prove the main theorem.
\begin{ttheorem}[5C]
Let $G$ be a hyperbolic group, and let $H\subset G$ be a quasiconvex subgroup.  If $[G:H] = \infty$, then there exists a (set-theoretic) section $s:G/H \to G$ of the quotient map $G \to G/H$ such that $s(G/H)$ is a net in $G.$
\end{ttheorem}
\begin{proof}
We first claim that there exists a constant $C_2$ such that for any $g\in G$, there exists a point $g'\in G$ and a geodesic ray $r$ from $1$ through $g'$ such that $d(g, g') < C_2$ and $d(r(t), H)$ is unbounded as $t\to\infty.$  By Proposition~\ref{proposition:3G}, there exists a point $g'$ with $d(g, g') < C$ for some constant $C$ (independent of $g$) and a geodesic ray $r$ from $1$ that passes through $g'.$  Suppose $r\subset \nbhd{l}{H}$ for some $l > 0$.  Set $x = r(\abs{g} + l)$, and choose some $h\in H$ such that $d(x, h) \leq l$.  Then
\begin{align*}
x.h
	&= \half\left(\abs{x} + \abs{h} - d(x, h)\right)
	\geq \abs{x} - d(x, h)
	\geq \abs{g}.
\end{align*}
Since $H$ is $K$-quasiconvex, there exists a geodesic $[1h]$ in $G$ lying in $\nbhd{K}{H}$.  Any two geodesics $[1g'], [1h]$ stay a distance at most $4\delta$ apart until time $g'.h$.  Hence 
\begin{align}\label{eq:close-to-H}
d(g, H) \leq d(g, [1h]) + K \leq d(g, g') + d(g', [1h]) + K + \leq K + C + 4\delta.
\end{align}
Set $C' = K + C + 4\delta$.  If $d(g, H) > C'$, contradicting \eqref{eq:close-to-H}, then the distance $d(r(t), H)$ must be unbounded; the claim then holds with $C_2 = C$.  Suppose instead that $d(g, H) \leq C'.$  Choose some $p\in \nbhd{C'}{1}$ such that $g \in Hp.$  Since $[G:H]$ is infinite, there exists some $t\in G$ with $d(t, H) > C'$.  Fix some such $t$ minimizing $\abs{t}$.  The point $gp^{-1}t$ then satisfies $d(gp^{-1}t, H) \geq d(Ht, H) \geq d(t, H) > C'$ and $d(gp^{-1}t, g) \leq C_2$, where
\begin{align*}
C_2 = C'+ \min\set{\abs{t'}:\, t'\in G, t'\not\in \nbhd{C'}{H}}.
\end{align*}
The claim therefore holds for all $g\in G$.

Fix $g\in G$.  By the claim above, there exist $g', x\in G$ such that $d(g, g') < C_2$, $g'\leq x$, and $d(x, H) \geq \abs{g} + \half C_1.$  Choose a point $x'\in S\cap xH$.  The Gromov product $x.x'$ satisfies
\begin{align*}
x.x' = \half\left(\abs{x} + \abs{x'} - d(x, x')\right)
	\geq \abs{x'} - \half C_1
	\geq \abs{g}
\end{align*}
by Lemma~\ref{lemma:5A}.  Any two geodesics $[1x]$ and $[1x']$ remain a distance no greater than $4\delta$ apart until time $x.x'$, so $d(g, [1x']) \leq d(g, g') + d(g', [1x']) < C_2 + 4\delta$. Set $L = \set{x\in \Lambda:\, \overline{x}\in S}$.  Then any $g\in G$ satisfies $d(g, \pi(\overline{L})) < C_2 + 4\delta$, where 
\begin{align*}
\overline{L} &= \set{x\in \Sigma^*\setvert \exists y\in \Sigma^*:\, xy\in L}
\end{align*}
denotes the prefix closure of $L$.  By Lemma~\ref{lemma:5I}, $\overline{L}$ is regular.  Let $M$ be a deterministic finite automaton with $L(M) = L$.  Since $M$ has only finitely many states, any word in $\overline{L}$ is within a bounded distance of a word in $L$; explicitly, any $x\in \overline{L}$ satisfies $d(x, L) < C_3$, where $C_3$ is the number of states of $M.$  Hence
\begin{align*}
d(g, L) \leq d(g, g') + d(g', L) < C_2 + C_3 + 4\delta
\end{align*}
By Lemma~\ref{lemma:5A}, each coset $g_0 H$ contains at most $\# \nbhd{C_1}{1}$ elements of $S$.  Choosing exactly one point in each intersection $S\cap g_0 H$ produces a section $s:G/H \to G$ such that $d(p, s(G/H)) < \#\nbhd{C_1}{1} + C_2 + C_3 + 4\delta$ for all vertices $p\in \Gamma$.
\end{proof}
Since the image $S = s(G/H)$ of the section $s:G/H \to G$ in Theorem~\ref{theorem:5C} is a net, it is also a hyperbolic metric space.  The left action of $G$ on the right coset space $G/H$ induces an action on $S$, given by $g.s(g') = s(gg')$.  By considering the corresponding homeomorphisms of the boundary induced by this action, we prove two results about the intersection of conjugate subgroups of $H$ below.  Set $H^g = g^{-1}Hg$ for any $g\in G$.  This conjugate depends only on the image of $g$ in the left coset space $H\action G$.  As such, we write $H^\gamma = \gamma^{-1} H \gamma = H^g$ for a coset $\gamma = Hg \in H\action G$.  In \cite{quasiconvex}, it is proved that any quasiconvex subgroup $H$ of a hyperbolic group $G$ has finite width; that is, $H^\gamma\cap H$ is finite for all but finitely many cosets $\gamma\in G/H$.  Using a completely different method, the section $s:G/H \to G$ of Theorem~\ref{theorem:5C}, we prove a weaker version of this result.  Specifically, we show in Proposition~\ref{proposition:5H} below that quasiconvex subgroups of infinite index in $G$ contain no infinite groups normal in $G$.  We require the following elementary lemma:
\begin{tlemma}[5G]
If $G$ is a finite extension of an infinite cyclic group, then any infinite cyclic subgroup $H\subset G$ has finite index.
\end{tlemma}
\begin{proof}
Choose an exact sequence $1 \to N \to G \xrightarrow{\pi} Q \to 1$ with $N$ cyclic and $Q$ finite.  Then $Q$ contains $\pi(H) = H/(H\cap N)$, so $[H:H\cap N]$ is finite.  In particular, $H\cap N$ is non-trivial.  Thus $[N:H\cap N]$ is finite.  The index $[G:H] \leq [G:H\cap N] = [G:N][N:H\cap N] = \#Q [N:H\cap N]$ is therefore also finite.
\end{proof}
We also need the following result, which is interesting independently of its use in proving Proposition~\ref{proposition:5H}.
\begin{tproposition}[5F]
For any $g\in G$, let $L_g$ denote the isometry $L_g(g') = gg'$ on the vertices of $\Gamma$.  Extend $L_g$ to a graph automorphism of $\Gamma$.  Suppose $G$ is not elementary.  Then the homomorphism $g \to (L_g)_\infty$ has finite kernel.
\end{tproposition}
\begin{proof}
Let $K$ denote the group of $g\in G$ with $(L_g)_\infty = \id$, and fix $g\in K$.  We first claim that there exists some constant $N = N(g)$ such that $[x^N, g] = 1$ for all $x\in G$.  By Proposition~\ref{proposition:1J}, the supremum $N_0$ of the order of all torsion elements of $G$ is finite.  The claim therefore holds immediately for all torsion $x\in G$ with $N = N_0!$.  Thus let $x\in G$ be an arbitrary element of infinite order.  By Proposition~\ref{proposition:3G}, there exists a constant $C$, independent of $g$ and $x$, such that $d(r, x) \leq C$ for some geodesic $r$.  Since $(L_g)_\infty$ acts trivially on $(r_t)\in \del G$, the distance $d(gr_t, r_t)$ is bounded.  Both $[1g]$ and $gr$ are geodesics, so the union $[1g] \cup gr$ is a $(1, 2\abs{g})$-quasigeodesic.  Thus there exists some geodesic ray $r'$ from $1$ such that $d(gr_t, r'_t) \leq C'$ for some constant $C' = C'(g)$ by Lemma~\ref{lemma:1F}.  The distance $d(r_t, r'_t) \leq d(r_t, gr_t) + d(gr_t, r'_t)$ is then bounded.  But 
\begin{align*}
r_t.r'_t = \half \left(\abs{r_t} + \abs{r'_t} - d(r_t, r'_t)\right) = t - d(r_t, r'_t),
\end{align*}
so $d(r_t, r'_t) \leq 4\delta$ for all time $t$.  Choosing some $t$ with $d(r_t, x) \leq C$, we therefore have
\begin{align*}
d(gx, x)
	&\leq d(gx, gr_t)+ d(gr_t, r_t) + d(r_t, x) \\
	&= 2d(x, r_t) + d(gr_t, r_t) \\
	&\leq 2d(x, r_t) + d(gr_t, r'_t) + d(r_t, r'_t) \\
	&\leq 2C + C' + 4\delta.
\end{align*}
Thus $\fabs{x^{-1}gx} \leq 2C + C' + 4\delta$ for all $x\in G$ of infinite order.  Fix such an $x\in G$.  Then $\abs{x^{-n}gx^n}\leq 2C+C'+4\delta$ for all $n > 0$, so there exist distinct $n, m$ such that $x^{-n}gx^n = x^{-m}gx^{m}$ and $0 < n < m < K$, where $K = \#\nbhd{2C+C'+4\delta}{1}$.  Hence any $x\in G$ of infinite order satisfies $[x^{K!}, g] = 1$.  The claim therefore holds for arbitrary $x\in G$ with $N = N_0! K!.$

By Proposition~\ref{proposition:1I}, there exist $x_1, x_2\in G$ such that $\fgroupgen{x_1, x_2}$ is a free group of rank $2$.  The commutators $[x_1^N, g]$ and $[x_2^N, g]$ vanish by the claim above.  The centralizers $C_G(x_1^N)$ and $C_G(x_2^N)$ are finite extensions of $\fgroupgen{x_1^N}$ and $\fgroupgen{x_2^N}$, respectively, by Proposition~\ref{proposition:1H} and Lemma~\ref{lemma:5G}.  Suppose $\groupgen{g}$ is infinite.  By Lemma~\ref{lemma:5G}, both $\fgroupgen{g}$ and $\fgroupgen{x_i^N}$ have finite index in $C_G(x_i^N)$ for each $i=1, 2$.  It follows that $\fgroupgen{g}\cap \fgroupgen{x_1^N}\cap \fgroupgen{x_2^N}$ has finite index in $\fgroupgen{g}$.  But $\fgroupgen{x_1^N, x_2^N}\subset \fgroupgen{x_1, x_2}$ is free, so $\fgroupgen{x_1^N} \cap \fgroupgen{x_2^N}$ is trivial.  It follows that $g$ must have finite order.

Thus $K$ consists entirely of torsion.  Let $g_1, \dots, g_r\in K$ be distinct representatives of the conjugacy classes of torsion in $G$ that intersect $K$.  By the claim above, there exist constants $N(g_i)$ for $i=1, \dots, r$ such that $[x^{N(g_i)}, g_i] = 1$ for all $x\in G$.  Set $N = N(g_1) \cdots N(g_r)$, and let $H = \groupgen{x^N:\, x\in G}$.  Then $H\normal G$ and $[H, g_i] = 1$ for all $g_i$.  Write $x^y$ for the conjugate $y^{-1}xy$.  Then for all $x\in G$ and $h\in H$,
\begin{align*}
(g_i^x)^h
	= h^{-1}g_i^x h
	= \left((h^{x^{-1}})^{-1} g_i (h^{x^{-1}}) \right)^x
	= g_i^x
\end{align*}
for each $g_i$.  Thus $[h, g_i^x] = 1$ for all $h\in H$ and $x\in G$.  The commutator $[h, k]$ therefore vanishes for all $k\in K$.  By Proposition~\ref{proposition:1I}, $G$ contains an element $x\in G$ of infinite order.  Then the centralizer $C_G(x^N)$ is a finite extension of $\Z$ by Proposition~\ref{proposition:1H}.  Hence $\fgroupgen{x^N}$ has finite index in $C_G(x^N)$ by Lemma~\ref{lemma:5G}.  But $C_G(x^N)$ contains $\fgroupgen{x^N}\cross K$, so
\begin{align*}
\#K 
	= \left[\groupgen{x^N} \cross K : \groupgen{x^N}\right] 
	\leq \left[C_G(x^N):\groupgen{x^N}\right] < \infty,
\end{align*}
as required.
\end{proof}
\begin{tproposition}[5H]
Let $G$ be a hyperbolic group, and let $H\subset G$ be a quasiconvex subgroup with $[G:H] = \infty$.  Let $K\subset H$ with normal closure $K^G$ in $G$.  If $K$ is infinite, then $[K^G:K] = \infty$.  In particular, any subgroup of $H$ normal in $G$ is finite.
\end{tproposition}
\begin{proof}
The corollary is trivial if $G$ is finite or a finite extension of $\Z$, so assume without loss of generality that $G$ is non-elementary.  Suppose instead that $K^G/K$ is a finite set of order $n$.  Each element of $G$ acts on the coset space $K^G/K$ by conjugation, giving a homomorphism $\rho:G \to S_n$.  The kernel $G' = \ker \rho$ has finite index in $G$, so $\del G' = \del G$.  Replacing $G$ by the quasi-isometric group $G'$, we can therefore assume that $n = 1$; that is, $K\normal G$.

Let $s:G/H \to G$ denote the section given by Theorem~\ref{theorem:5C}.  Then $G$ acts on $S = s(G/H)\subset G$ by $g.s(x) = s(gx)$.  Since $s$ is an injective function from the right coset space $G/H$ to $G$, this action is well-defined.  The stabilizer of any $s(x) \in S$ is $H^x \supset K$, so $K$ acts trivially on $S$.  Since $S$ is a net in $G$, the inclusion $i:S\to G$ is a $(1, \epsilon)$-quasi-isometry for some $\epsilon$.  Hence $i$ induces a bijection (in fact, a homeomorphism) $i_\infty:\del S \to \del G$.  Since $s(x)$ minimizes $\abs{\cdot}$ in $x H$, we have
\begin{align*}
\fabs{s(x)} \leq \fabs{g^{-1}s(gx)} \leq \fabs{g^{-1}} + \fabs{s(gx)} = \fabs{g} + \fabs{s(gx)}.
\end{align*}
Thus $\abs{s(gx)} \geq \abs{s(x)} - \abs{g}$.  By Lemma~\ref{lemma:5A},
\begin{align}
d(g.s(x), gs(x))
	&= d(s(gx), gs(x)) \nonumber\\
	&\leq C_1 + \abs{sg(x)} + \abs{gs(x)} - 2\abs{s(gx)} \nonumber\\
	&\leq C_1 + 2\abs{g} \label{eq:action}
\end{align}
for all $x$.  Thus the diagram
\begin{align}\label{eq:cd}
\xymatrix@C+.75cm@R+.375cm{
	G \ar[r]^(.35){g\to (L'_g)_\infty} \ar[dr]_(.4){g\to (L_g)_\infty} & \Homeo(\del S) \\
	& \Homeo(\del G) \ar[u]_(.46){(i_\infty)^*}
}
\end{align}
commutes, where $L_g(g') = gg'$ and $L'_g s(x) = g.s(x)$.  We hence have $(L_g)_\infty = 1$ for all $g\in K$.  By Proposition~\ref{proposition:5F}, $K$ is finite.
\end{proof}
\bibliography{References}
\end{document}